\documentclass[11pt,a4paper]{article}
\usepackage[english]{babel}
\usepackage[english]{layout}
\usepackage{amsmath,bbm,latexsym}
\usepackage{amsfonts,amssymb}
\usepackage{color,graphics,graphicx,amsthm}
\usepackage[T1]{fontenc}
\usepackage{color}
\def\eb{\begin{eqnarray*}}
\def\ee{\end{eqnarray*}}
\newcommand{\bc}{\begin{center}}
\newcommand{\ec}{\end{center}}

\renewcommand{\l}{\lambda}
\newtheorem{thm}{Theorem}[section]
\newtheorem{defi}{Definition}[section]

\def\non{\nonumber}

\def\mf{\mathfrak}
\newcommand{\al}{\alpha}
\newcommand{\be}{\beta}
\newcommand{\ga}{\gamma}
\newcommand{\de}{\delta}
\newcommand{\ep}{\varepsilon}
\newcommand{\et}{\eta}
\newcommand{\tf}{\vartheta}

\begin{document}
\title{Creation operators and algebraic Bethe ansatz for the  elliptic quantum group $E_{\tau,\eta}(so_3)$}
\author{Nenand Manojlovi\'c\footnote{e-mail address: nmanoj@ualg.pt} and Zolt\'an Nagy\footnote{e-mail address: znagy@ualg.pt} \vspace{5pt}\\ Departamento de Matem\'atica \\ FCT,  Campus de Gambelas\\
Universidade do Algarve \\8005-139 Faro, Portugal}
\date{}
\maketitle
\abstract{We define the elliptic quantum group $E_{\tau,\eta}(so_3)$ and the transfer matrix
corresponding to its simplest highest weight representation. We use Bethe anstaz method to construct the
creation operators
as polynomials of the Lax matrix elements expressed through a recurrence relation. We give common eigenvectors and
eigenvalues of the family of commuting transfer matrices.}
\section{Introduction}
In this article, we report new results on the application of \emph{algebraic} Bethe ansatz to
the elliptic (or dynamical) quantum group $E_{\tau,\eta}(so_3)$. The elliptic quantum group is
the algebraic structure associated to elliptic solutions of the star-triangle relation.
This equation appears in interaction-round-a-face models in statistical mechanics. As it was shown
by Felder \cite{Fe}, this structure is also related to Knizhnik-Zamolodchikov-Bernard equation of conformal
field theory on tori. Moreover, in a different direction, to each solution of the (see \cite{Ji}) star-triangle relation
a dynamical $R$-matrix can be associated. This $R$-matrix, in turn, will define
an algebra similar to quantum groups appearing in the quantum inverse scattering method (QISM), it is actually
a quasi-Hopf deformation of the more familiar quantum group structure \cite{Ji2}.
Despite all the differences,
this new structure preserves a prominent feature of quantum groups: a tensor product of representations
can be defined.

The adjective dynamical refers to the fact that the $R$-matrix appearing in these
structures contains a parameter which in the classical limit will be interpreted as the position coordinate
on the phase space of a classical system and the resulting classical $r$-matrix will depend on it. In
the quantum setting, apart from the appearance of this extra parameter the Yang-Baxter equation (YBE) is
also deformed. At the technical level, the main difference between usual quantum groups and the one
we are about to describe lies not so much in the elliptic nature of the appearing functions as rather in
the introduction of the extra "dynamical" parameter and the corresponding deformation of the YBE.

In QISM, the physically interesting quantity is the transfer matrix. The hamiltonian of the
model and other observables are derived from it. The knowledge of its spectrum is thus essential.
Different kinds of methods under the federating name of Bethe
ansatz have been developed to calculate the eigenvalues of the transfer matrix \cite{Fa,Ko,Ku}.
The question whether the algebraic Bethe ansatz (ABA) technique can be applied to transfer matrices
appearing in the context of dynamical quantum groups has received an affirmative answer from Felder and
Varchenko \cite{Fe2,FeVa}.
They showed how to implement ABA for the elliptic quantum group $E_{\tau,\eta}(sl_2)$, they also showed its
applications to IRF models and Lam\'e equation. Later, for the $E_{\tau,\eta}(sl_n)$ elliptic
quantum group the nested Bethe ansatz method was used \cite{Es,Sa}
and a relation to Ruijsenaars-Schneider \cite{Sa} and quantum Calogero-Moser Hamiltonians was established \cite{ABB}.

In the first section of this paper we introduce the basic definitions of dynamical $R$-matrix,
Yang-Baxter equation, representations,
operator algebra and commuting transfer matrices.
We define elements $\Phi_n$ in the operator algebra which have the necessary symmetry properties to be the
creation operators of the corresponding Bethe states. As it turns out, the creation operators are not
simple functions of the Lax matrix entries, unlike in \cite{FeVa}, but they are complicated
polynomials of three generators $A_1(u), B_1(u), B_2(u)$ in the elliptic operator algebra. We give the
recurrence relation which  defines the creation operators. Moreover, we fully implement the algebraic Bethe ansatz
on the simplest example of highest weight module.
 That is we calculate the action of the transfer matrix on the Bethe vectors and from the vanishing
of the unwanted terms we derive the Bethe equations. We also give the
explicit formulas for the corresponding eigenvalues.

\section{Representations of $E_{\tau,\eta}(so_3)$ and transfer matrices}
\subsection{Definitions}
Let us first recall the basic definitions which will enter our construction. First, we fix two complex numbers
$\tau,\eta$ such that $Im(\tau) > 0$.
The central object in this paper is the $R$-matrix $R(q,u)$ which depends on two arguments $q,u \in \mathbb{C}$:
the first one is referred to
as the dynamical parameter, the second one is called the spectral parameter. The elements of the $R$-matrix
are written in terms of Jacobi's theta function:
\eb
\tf(u)=-\sum_{j\in \mathbb{Z}}\exp \left(\pi i\left(j+\frac{1}{2}\right)^2+2\pi i \left(j+\frac{1}{2}\right)
\left(u+\frac{1}{2}\right)\right)
\ee
This function has two essential properties. It is quasiperiodic:
\eb
\tf(u+1)=-\tf(u); \qquad \tf(u+\tau)=-e^{-i \tau-2 i u}\tf(u)
\ee
and it verifies the identity:
\eb
&&\tf(u+x)\tf(u-x)\tf(v+y)\tf(v-y)=\tf(u+y)\tf(u-y)\tf(v+x)\tf(v-x)+\\
&&\tf(u+v)\tf(u-v)\tf(x+y)\tf(x-y)
\ee

The entries of the $R$-matrix are written in terms of the following functions.

\begin{eqnarray*}
g(u)  &=&  \frac{\tf(u-\et) \tf(u-2\et)}{\tf(\et) \tf(2 \et)}\\
\al(q_1,q_2,u) & =&  \frac{\tf(\et-u)\tf(q_{12}-u)}{\tf(\et)\tf(q_{12})}\\
\be(q_1,q_2,u) &=&  \frac{\tf(\et-u)\tf(u)\tf(q_{12}-2\et)}{\tf(-2\et)\tf(\et)\tf(q_{12})}\\
\ep(q,u)&=& \frac{\tf(\et+u)\tf(2\et-u)}{\tf(\et)\tf(2\et)}-\frac{\tf(u)\tf(\et-u)}{\tf(\et)\tf(2\et)}
\left( \frac{\tf(q+\et)\tf(q-2\et)}{\tf(q-\et)\tf(q)}+\frac{\tf(q-\et)\tf(q+2\et)}{\tf(q+\et)\tf(q)}\right)\\
\ga(q_1,q_2,u)&= & \frac{\tf(u)\tf(q_1+q_2-\et-u)\tf(q_1-2\et)\tf(q_2+\et)}{\tf(\et)\tf(q_1+q_2-2\et)\tf(q_1+\et)\tf(q_2)}\\
\de(q,u) & = & \frac{\tf(u-q)\tf(u-q+\et)}{\tf(q)\tf(q-\et)}
\end{eqnarray*}

The $R$-matrix itself will act on the tensor product $V \otimes V$ where $V$ is a three-dimensional
complex vector space with the standard basis $\{e_1,e_2,e_3\}$. The matrix units $E_{ij}$ are defined
in the usual way: $E_{ij}e_k=\delta_{jk}e_i$. We will also need the following diagonal matrix later on: $h=E_{11}-E_{33}$.

Now we are ready to write the explicit form of the $R$-matrix.
\eb\label{Rmat}
R(q,u)&=&g(u)E_{11}\otimes E_{11}+g(u)E_{33}\otimes E_{33}+\ep(q,u)E_{22}\otimes E_{22}\\
&+&\al(\eta,q,u)E_{12}\otimes E_{21}+\al(q,\eta,u)E_{21}\otimes E_{12}+\al(-q,\eta,u)E_{23}\otimes E_{32}\\
&+& \al(\eta,-q,u)E_{32}\otimes E_{23}\\
&+& \be(\eta,q,u) E_{11}\otimes E_{22}+\be(q,\eta,u) E_{22}\otimes E_{11}+\be(-q,\eta,u) E_{22}\otimes E_{33}\\
&+& \be(\eta,-q,u)E_{33}\otimes E_{22}\\
&+&\ga(-q,q,u)E_{11}\otimes E_{33}+\ga(-q,\eta,u)E_{12} \otimes E_{32}- \ga(\eta,q,u) E_{21} \otimes E_{23}\\
&+& \ga(q,-q,u) E_{33} \otimes E_{11}+ \ga(q,\eta,u) E_{32} \otimes E_{12}- \ga(\eta,-q,u) E_{23} \otimes E_{21}\\
&+& \de(q,u) E_{31}\otimes E_{13}+\de(-q,u) E_{13} \otimes E_{31}
\ee
This $R$-matrix also enjoys the unitarity property:
\begin{equation}\label{unit}
R_{12}(q,u)R_{21}(q,-u)=g(u)g(-u)\mathbbm{1}
\end{equation}
and it is of zero weight:
\eb
\left[h \otimes \mathbbm{1}+\mathbbm{1} \otimes h,R_{12}(q,u)\right]=0 \qquad  (h \in \mf{h})
\ee

The $R$-matrix also obeys the dynamical quantum Yang-Baxter equation (DYBE) in $End(V \otimes V \otimes V)$:
\eb
&&R_{12}(q-2\eta h_3,u_{12})  R_{13}(q,u_1) R_{23}(q-2\eta h_1,u_2)=\\
 && R_{23}(q,u_2)R_{13}(q-2\eta h_2,u_1)
R_{12}(q,u_{12})
\ee
where the "dynamical shift" notation has the usual meaning:
\begin{eqnarray}\label{shift}
R_{12}(q-2\eta h_3,u) \cdot v_1\otimes v_2 \otimes v_3 = \left(R_{12}(q-2\eta \lambda,u) v_1\otimes v_2\right)
\otimes v_3
\end{eqnarray}
whenever $h v_3= \lambda v_3$. Shifts on other spaces are defined in an analogous manner.

Let us also describe a more intuitive way of looking at this shift. Define first
the shift operator acting on functions of the
dynamical parameter:
\eb
\exp(2\eta \partial_q) f(q)=f(q+2\eta) \exp(2\eta \partial_q)
\ee
Then equation (\ref{shift}) can also be written in the following form:
\eb
R_{12}(q-2\eta h_3,u) = \exp(-2\eta h_3 \partial_q) R_{12}(q,u) \exp(2\eta h_3 \partial_q)
\ee
in the sequel we will use whichever definition is the fittest for the particular point in our calculation.

\subsection{Representation; operator algebra}
Now we describe the notion of representation of (or module over) $E_{\tau,\eta}(so_3)$. It
 is a pair $(\mathcal{L}(q,u),W)$ where $W$ is a diagonalizable $\mf{h}$-module, that is, $W$ is a direct sum of
 the weight subspaces
$W=\oplus_{\l \in \mathbb{C}}W[\l]$ and $\mathcal{L}(q,u)$ is an operator in $\mathrm{End}(V \otimes W)$ obeying:
\eb\label{RLL}
&&R_{12}(q-2\eta h_3,u_{12})  \mathcal{L}_{13}(q,u_1) \mathcal{L}_{23}(q-2\eta h_1,u_2)=\\
 && \mathcal{L}_{23}(q,u_2)\mathcal{L}_{13}(q-2\eta h_2,u_1)
R_{12}(q,u_{12})
\ee

$\mathcal{L}(q,u)$ is also of zero weight
\eb
\left[h_V \otimes \mathbbm{1}+\mathbbm{1} \otimes h_W , \mathcal{L}_{V,W}(q,u)\right]=0 \qquad (h \in \mf{h})
\ee
where the subscripts remind the careful reader that in this formula $h$ might act in a different way on spaces
$W$ and $V$.

An example is given immediately by $W=V$ and $\mathcal{L}(q,u)=R(q,u-z)$ which is called the fundamental
representation with evaluation point $z$ and is denoted by $V(z)$.
A tensor product of representations can also be defined which corresponds to the existence of a coproduct-like
structure at the abstract algebraic level. Let $(\mathcal{L}(q,u),X)$ and $(\mathcal{L}'(q,u),Y)$
be two $E_{\tau,\eta}(so_3)$
modules, then $\mathcal{L}_{1X}(q-2\eta h_Y,u)\mathcal{L}_{1Y}(q,u),X\otimes Y$ is a
representation of $E_{\tau,\eta}(so_3)$ on
$X \otimes Y$ endowed, of course, with the tensor product $\mf{h}$-module structure.

The operator $\mathcal{L}$ is reminiscent of the quantum Lax matrix in the FRT formulation
of the quantum inverse scattering
method, although it obeys a different exchange relation,
therefore we will also call it a Lax matrix. This allows us to view
the $\mathcal{L}$ as a matrix with operator-valued entries.

Inspired by that interpretation, for any module over $E_{\tau,\eta}(so_3)$ we define the corresponding
\textbf{operator algebra} of finite difference operators.
Let us take an arbitrary representation $\mathcal{L}(q,u) \in \mathrm{End}(V \otimes W)$.
The elements of the operator algebra corresponding to this representation will act on the space $\mathrm{Fun}(W)$ of
meromorphic functions of $q$ with values in $W$. Namely let $L \in \mathrm{End}(V \otimes \mathrm{Fun}(W))$
be the operator defined as:
\begin{equation}\label{Lti}
L(u)=\left( \begin{array}{ccc}
A_1(u)& B_1(u)& B_2(u)\\
C_1(u) & A_2(u) & B_3(u)\\
C_2(u) & C_3(u) &A_3(u)
 \end{array}\right)=\mathcal{L}(q,u)e^{-2\eta h \partial_q}
\end{equation}
We can view it as a matrix with entries in $\mathrm{End}(\mathrm{Fun}(W))$:
It follows from equation (\ref{RLL}) that $L$
verifies:
\begin{eqnarray}\label{RLLti}
R_{12}(q-2\eta h,u_{12}) \ L_{1W}(q,u_1) L_{2W}(q,u_2)= L_{2W}(q,u_2)L_{1W}(q,u_1) \
\tilde{R}_{12}(q,u_{12})
\end{eqnarray}
with $\tilde{R}_{12}(q,u):= \exp(2\eta(h_1+h_2)\partial_q)R_{12}(q,u)\exp(-2\eta(h_1+h_2)\partial_q)$

The zero weight condition on $L$ yields the relations:
\eb\label{order}
\left[h,A_i\right]=0 ; \ \ \left[h,B_j\right]=-B_j \quad (j=1,3), \ \left[h,B_2\right]=-2B_2\\
\left[h,C_j\right]=C_j \quad (j=1,3), \ \left[h,C_2\right]=2C_2
\ee
so $B_i$'s act as lowering and $C_i$'s as raising operators.
From the definition \eqref{Lti} one can derive the action of the operator algebra generators on functions:
\eb
A_1(u) f(q)= f(q-2\eta)A_1(u);\ B_1(u)f(q)= f(q)B_1(u); \\
B_2(u)f(q)= f(q+2\eta)B_2(u)
\ee
and analogously for the other generators.
We display here those commutation relations which are necessary for the construction  of the Bethe
vectors, the remaining ones can be extracted from \eqref{RLLti}.
\eb
B_1(u_1)B_1(u_2)&=&\omega_{21}\left(B_1(u_2)B_1(u_1)-\frac{1}{y_{21}(q)}B_2(u_2)A_1(u_1)\right)+
\frac{1}{y_{12}(q)}B_2(u_1)A_1(u_2) \label{crB1B1}  \\
A_1(u_1)B_1(u_2)&=&z_{21}(q)B_1(u_2)A_1(u_1)-\frac{\al_{21}(\eta,q)}{\be_{21}(q,\eta)}B_1(u_1)A_1(u_2)  \\
A_1(u_1)B_2(u_2)&=&\frac{1}{\ga_{21}(q,-q)}\left( g_{21}B_2(u_2)A_1(u_2)+\ga_{21}(\eta,-q)B_1(u_1)B_1(u_2) \non
-\de_{21}(-q)B_2(u_1)A_1(u_1)\right) \\
 B_1(u_2)B_2(u_1)&=&\frac{1}{g_{21}}\left( \be_{21}(-q,\eta)B_2(u_1)B_1(u_2)+\al_{21}(\eta,-q)B_1(u_1)B_2(u_2)\right)
 \\
 B_2(u_2)B_1(u_1)&=&\frac{1}{g_{21}}\left( \be_{21}(\eta,-q)B_1(u_1)B_2(u_2)+\al_{21}(-q,\eta)B_2(u_1)B_1(u_2)\right)\label{crB2B1}
\ee
where
\eb
\omega(q,u)&=&\frac{g(u)\ga(q,-q,u)}{\ep(q,u) \ga(q,-q,u)+\ga(q,\eta,u)\ga(\eta,-q,u)}\\
y(q,u)&=&\frac{\ga(q,-q,u)}{\ga(q,\eta,u)} \nonumber\\
z(q,u)&=&\frac{g(u)}{\be(q,\eta,u)}\nonumber
\ee
and as usual
\[
y_{12}(q)=y(q,u_1-u_2) \ \ \textrm{etc.}
\]
Furthermore, the function $\omega(q,u)$ is actually independent of $q$, a property which will prove important later on,
and takes the following simple form:
\eb\label{omeg}
\omega(u)=\frac{\tf(u+\eta)\tf(u-2\eta)}{\tf(u-\eta)\tf(u+2\eta)}
\ee
This function also verifies the following property:
\eb
\omega(u)\omega(-u)=1 \ .
\ee
Finally the following theorem shows how to associate a family of commuting quantities to a representation
of the elliptic quantum group
\begin{thm}
Let $W$ be a representation of $E_{\tau,\eta}(so_3)$. Then the transfer matrix defined by $t(u)=Tr \tilde{L}(u) \in
\mathrm{End}(\mathrm{Fun}(W))$
preserves the subspace $\mathrm{Fun}(W)[0]$ of functions with values in the zero weight subspace of $W$.
When restricted to this subspace, they commute at different values of the spectral parameter:
\eb
\left[t(u),t(v)\right]=0
\ee
\end{thm}
\begin{proof}
The proof is analogous to references \cite{FeVa3,ABB}
\end{proof}
\section{Bethe ansatz}
In this section we fix a highest weight module $W=V(z_1)\otimes \ldots \otimes V(z_n)$.
The vector $|0\rangle =e_1\otimes \ldots \otimes e_1 \in \textrm{Fun}(W)$ is
 a highest weight vector of weight $n$ of this module, and every highest weight vector can be written in the
 form $|\Omega \rangle = f(q) |0 \rangle$ with a non-zero meromorphic function $f$. We have indeed:
 \eb
 C_i(u)|\Omega \rangle=0  \qquad (i=1,2,3)
 \ee
showing that $|\Omega \rangle$ is a highest weight vector.

\eb
A_1(u)|\Omega\rangle=a_1(u)\frac{f(q-2\eta)}{f(q)}|\Omega\rangle\\
\quad A_2(u)|\Omega\rangle=a_2(q,u)|\Omega\rangle \quad A_3(u)|\Omega\rangle=a_3(q,u)\frac{f(q+2\eta)}{f(q)}|\Omega\rangle
\ee
with the eigenvalues:
\eb
a_1(u)&=&\prod_{i=1}^n\frac{\tf(u-z_i-\et) \tf(u-z_i-2\et)}{\tf(\et) \tf(2 \et)}\\
a_2(q,u)
&=&\frac{\tf(q-2\eta n-\eta)}{\tf(q-\eta)}\prod_{i=1}^n\frac{\tf(u-z_i-\eta)\tf(u-z_i)}{\tf(\et)\tf(2\et)}\\
a_3(q,u)&=&\frac{\tf(q-2\eta n)\tf(q-2\eta n+\eta)}{\tf(q+\eta)\tf(q)}\prod_{i=1}^n
\frac{\tf(u-z_i)\tf(u-z_i+\eta)}{\tf(\eta)\tf(2\eta)}
\ee
Notice that $a_1$ is independent of $q$.
It is easy to see, that the
zero weight subspace $W[0]$ is nontrivial for this module.

In this situation we cannot look for the common eigenvectors of $t(u)$ in the form $\Phi(u_1,\ldots,u_n)=
B_1(u_1)\ldots B_1(u_n)|\Omega\rangle$ since $B_1(u)B_1(v)\neq B_1(v)B_1(u)$ and the resulting Bethe vector
would not be symmetric under the interchange of the parameters $u_i$.

Instead, we should be inspired by Tarasov's implementation of the algebraic
Bethe ansatz to the Izergin-Korepin model \cite{Ta}. In that case
the $R$-matrix is nondynamical, but
has nonzero entries at the same positions as the dynamical $R$-matrix considered here. The comparison
suggests that the Bethe operator will contain $B_1(u_1)\ldots B_1(u_n)$ with coefficient $1$ but will also
have a "correction" to that expressed in terms of $B_2(u)$ and $A_1(u)$. It also suggests that the symmetry
under the interchange of spectral parameters is replaced by the property:
\eb
\Phi_n(u_1,\ldots,u_n)=\zeta_{i+1,i}\Phi_n(u_1, \ldots, u_{i-1},u_{i+1},u_{i},u_{i+2},\ldots,u_n) \qquad (i=1,\ldots,n-1)
\ee
with a function $\zeta$ to be determined. In the sequel we give the Bethe creation operator in a recurrence form
and describe its generalized symmetry property.

\begin{defi}
Let $\Phi_n$ be defined be the recurrence relation for $n\geq 0$:
\eb
&&\Phi_n(u_1,\ldots,u_n)=B_1(u_1)\Phi_{n-1}(u_2,\ldots, u_n)\\
&&-\sum_{j=2}^n\frac{\prod_{k=2}^{j-1}\omega_{jk}}{y_{1j}(q)}
\prod_{\substack{k=2\\k\neq j}}^n z_{kj}(q+2\eta)\ B_2(u_1) \Phi_{n-2}(u_2,\ldots,\widehat{u_j},\ldots,u_n)A_1(u_j)
\ee
where $\Phi_0=1; \ \Phi_1(u_1)=B_1(u_1)$ and the parameter under the hat is omitted.
\end{defi}

It may be useful to give explicitly the first three creation operators.
\eb
&&\Phi_1(u_1)=B_1(u_1)\\
&&\Phi_2(u_1,u_2)=B_1(u_1)B_1(u_2)-\frac{1}{y_{12}(q)}B_2(u_1)A_1(u_2)\\
&&\Phi_3(u_1,u_2,u_3)=B_1(u_1)B_1(u_2)B_1(u_3)-\frac{1}{y_{23}(q)}B_1(u_1)B_2(u_2)A_1(u_3)\\
&&-\frac{z_{32}(q+2\eta)}{y_{12}(q)}B_2(u_1)B_(u_3)A_1(u_2)-\frac{\omega_{32}z_{23}(q+2\eta)}{y_{13}(q)}B_2(u_1)
B_1(u_2)A_1(u_3)
\ee

The Bethe vector is then not completely symmetric under the interchange of two neighboring spectral parameters
but verifies the following property instead:
\eb
&&\Phi_2(u_1,u_2)=\omega_{21}\Phi_2(u_2,u_1)\\
&&\Phi_3(u_1,u_2,u_3)=\omega_{21}\Phi_3(u_2,u_1,u_3)=\omega_{32}\Phi_3(u_1,u_3,u_2)
\ee

For general $n$ we prove the following theorem.
\begin{thm}
$\Phi_n$ verifies the following symmetry property:
\begin{gather}\label{symm}
\Phi_n(u_1,\ldots,u_n)=\omega_{i+1,i}\Phi_n(u_1,\ldots,u_{i-1},u_{i+1},u_i,u_{i+2},\ldots,u_n)
\qquad (i=1,2,\ldots,n-1).
\end{gather}
\end{thm}
\begin{proof}
For the proof we refer to \cite{MaNa}.
%
\end{proof}

The next step in the application of the Bethe ansatz scheme is the calculation of the action of the transfer
matrix on the Bethe vector.  This will yield
three kinds of terms. The first part (usually called wanted terms in the literature)
will tell us the eigenvalue of the transfer matrix, the second part (called unwanted terms) must be annihilated
by a careful choice of the spectral parameters $u_i$ in $\Phi_n(u_1,\ldots,u_n)$; the vanishing of these unwanted
terms is ensured if the $u_i$ are solutions to the so called Bethe equations. The third part contains terms
ending with a raising operator acting on the pseudovacuum and thus vanishes.

The action of $A_1(u)$ on $\Phi_n$ is given by
\begin{eqnarray}\label{A1phi}
A_1(u)\Phi_n &=& \prod_{k=1}^n z_{ku}(q) \Phi_n A_1(u)+\\
&&\sum_{j=1}^n D_j \prod_{k=1}^{j-1}\omega_{jk} B_1(u)\Phi_{n-1}(u_1,\hat{u_j},u_n)A_1(u_j)+ \non \\
&&\sum_{l<j}^n E_{lj}\prod_{k=1}^{l-1} \omega_{lk} \non
\prod_{\substack{k=1\\k\neq l}}^{j-1}\omega_{jk} B_2(u) \Phi_{n-2}(u_1,\hat{u_l},\hat{u_j},u_n)A_1(u_l)A_1(u_j)
\end{eqnarray}
To calculate the first coefficients we expand $\Phi_n$ with the help
of the recurrence relation, then use the commutation relations to push $A_1(u_1)$ to the right.
This yields:
\eb
D_1&=&-\frac{\al_{1u}(\eta,q)}{\be_{1u}(q,\eta)}\prod_{k=2}^n
z_{k1}(q)\\
E_{12}&=&\left( \frac{\de_{1u}(-q)}{\ga_{1u}(q,-q)y_{12}(q-2\eta)}+
\frac{z_{1u}(q)\al_{2u}(\eta,q)\omega_{u 1}}{\be_{2u}(q,\eta)y_{u 1}(q)}\right)
\prod_{k=3}^n z_{k1}(q+2\eta)z_{k2}(q)
\ee
The direct calculation of the remaining coefficients is less straightforward. However, the symmetry of the left
hand side of \eqref{A1phi} implies that $D_j$ for $j\geq 1$ can be obtained by substitution
$u_1 \rightsquigarrow u_j$ in $D_1$ and  $E_{lj}$ by the substitution
$u_1 \rightsquigarrow u_l$, $ u_2 \rightsquigarrow u_j$

The action of $A_2(u)$ and $A_3(u)$ on $\Phi_n$ will yield also terms ending in $C_i(u)$'s.

The action of $A_2(u)$ on $\Phi_n$ will have the following structure.
\eb
A_2(u)\Phi_n  &=& \prod_{k=1}^n \frac{z_{u k}(q-2\eta(k-1))}{\omega_{u k}} \Phi_n A_2(u)+\\
&&\sum_{j=1}^n F^{(1)}_j \prod_{k=1}^{j-1}\omega_{jk} B_1(u)\Phi_{n-1}(u_1,\hat{u_j},u_n)A_2(u_j)+\\
&&\sum_{j=1}^n F^{(2)}_j \prod_{k=1}^{j-1}\omega_{jk} B_3(u)\Phi_{n-1}(u_1,\hat{u_j},u_n)A_1(u_j)+\\
&&\sum_{l<j}^n G^{(1)}_{lj} \prod_{k=1}^{l-1}\omega_{lk}\prod_{\substack{k=1\\k\neq l}}^{j-1} \omega_{jk}
B_2(u) \Phi_{n-2}(u_1,\hat{u_l},\hat{u_j},u_n) A_1(u_l)A_2(u_j)+\\
&&\sum_{l<j}^n G^{(2)}_{lj} \prod_{k=1}^{l-1}\omega_{lk}\prod_{\substack{k=1\\k\neq l}}^{j-1} \omega_{jk}
B_2(u)\Phi_{n-2}(u_1,\hat{u_l},\hat{u_j},u_n)A_1(u_j)A_2(u_l)+\\
&&\sum_{l<j}^n G^{(3)}_{lj}\prod_{k=1}^{l-1}\omega_{lk}\prod_{\substack{k=1\\k\neq l}}^{j-1} \omega_{jk}
B_2(u)\Phi_{n-2}(u_1,\hat{u_l},\hat{u_j},u_n)A_2(u_l)A_1(u_j)+\\
&&\textit{terms ending in C}
\ee

We give the coefficients $D^{(k)}_1$ and $E^{(k)}_{12}$, the remaining ones are obtained by the same
substitution as for $A_1(u)$

\eb
F^{(1)}_1&=&-{\frac{\al_{u 1}(q-2\eta,\eta)}{\be_{u 1}(q,\eta)}\prod_{k=2}^n
\frac{z_{1k}(q-2\eta(k-1))}{\omega_{1k}}}\\
F^{(2)}_1&=&\frac{1}{y_{u 1}(q)}\prod_{k=2}^n z_{k1}(q+2\eta)\\
G^{(1)}_{12}&=& \frac{1}{y_{u 1}(q)}\left( \frac{z_{u 1}(q)\al_{u 2}(q-2\eta,\eta)}{\be_{u 2}(q-2\eta,\eta)}-
\frac{\al_{u 1}(q,\eta)\al_{12}(q-2\eta,\eta)}{\be_{u 1}(q,\eta)\be_{12}(q-2\eta,\eta)}\right)
\prod_{k=3}^n \frac{z_{k1}(q+2\eta)z_{2k}(q-2\eta(k-1))}{\omega_{2k}}\\
G^{(2)}_{12}&=& \frac{\al_{u 1}(q,\eta)\al_{12}(q-2\eta,\eta)}{\be_{u 1}(q,\eta) y_{u 1}(q)\be_{12}(q-2\eta,\eta)}
\prod_{k=3}^n \frac{z_{k2}(q+2\eta)z_{1k}(q-2\eta(k-1))}{\omega_{1k}}\\
G^{(3)}_{12}&=&\frac{\al_{u 1}(q,\eta)}{\be_{u 1}(\eta,-q)}\left( \frac{z_{u 1}(q)}{\omega_{u 1}y_{u 2}(q)}-
\frac{\al_{u 1}(\eta,-q)}{y_{12}(q)\be_{u 1}(q,\eta)}\right)
\prod_{k=3}^n \frac{z_{k2}(q+2\eta)z_{1k}(q-2\eta(k-2))}{\omega_{1k}}
\ee
It is instructing to give explicitly the expression of $F^{(1)}_l$
\eb
F^{(1)}_l=-\frac{\al_{ul}(q-2\eta,\eta)}{\be_{ul}(q,\eta)}\times \frac{\tf(q-3\eta)}{\tf(q-2\eta n-\eta)}
\prod_{\substack{k=1\\k\neq l}}^{n} \frac{\tf(u_{lk}-2\eta)}{\tf(u_{lk}) \omega_{lk}}
\ee
The action of $A_3(u)$ on the Bethe vector is somewhat simpler.
\eb
A_3(u) \Phi_n &=& \prod_{k=1}^n \frac{\be_{u k}(\eta,-q)}{\ga_{u k}(q-2\eta(k-1),-)}
\Phi_n A_3(u)+ \\
&&\sum_{j=1}^n H_j \prod_{k=1}^{j-1}\omega_{jk} B_3(u)\Phi_{n-1}(u_1,\hat{u_j},u_n)A_2(u_j)+\\
&&\sum_{l<j}^n I_{lj}\prod_{k=1}^{l-1}\omega_{lk}\prod_{\substack{k=1\\k\neq l}}^{j-1} \omega_{jk}
B_2(u) \Phi_{n-2}(u_1,\hat{u_l},\hat{u_j},u_n) A_2(u_l)A_2(u_j)+\\
&&\textit{terms ending in C}
\ee
where to save space used the notation $\ga_{uk}(x,-)=\ga_{uk}(x,-x)$.
We give the coefficients $H_1$ and $I_{12}$, the rest can be obtained by the substitution of
the spectral parameters as before.
\eb
H_1&=&-\frac{1}{y_{u 1}(q)} \prod_{k=2}
\frac{z_{1k}(q-2\eta(k-2))}{\omega_{1k}}\\
I_{12}&=&\frac{1}{\ga_{u 2}(q,-q)}\left( \frac{\de_{u 2}(q)}{y_{12}(q-2\eta)}-\frac{\al_{u 1}(q,\eta)}{y_{u 2}(q-2\eta)}\right)
 \prod_{k=3}\frac{z_{2u}(q-2\eta (k-2))z_{1u}(q-2\eta (k-2))}{\omega_{1k}\omega_{2k}}
\ee

The next step is to find conditions for the cancelation of the unwanted terms.
We write the action of the transfer matrix in the following form:
\eb
&&t(u)\Phi_n |\Omega \rangle = \Lambda \Phi_n |\Omega\rangle +\\
&& \sum_{j=1}^n K^{(1)}_j \prod_{k=1}^{j-1}\omega_{jk} B_1(u) \Phi_n(u_1,\hat{u_j},u_n)+\\
&& \sum_{l<j}^n K^{(2)}_{lj} \prod_{k=1}^{l-1}\omega_{lk}\prod_{\substack{k=1\\k\neq l}}^{j-1} \omega_{jk}
B_2(u) \Phi_n(u_1,\hat{u_l},\hat{u_j},u_n)+\\
&& \sum_{j=1}^n K^{(3)}_j \prod_{k=1}^{j-1}\omega_{jk} B_3(u) \Phi_n(u_1,\hat{u_j},u_n)
\ee
The general form of the eigenvalue is written as
\eb
\Lambda(u,\{u_j\})&=&\prod_{k=1}^n z_{ku}(q)\times a_1(q,u)\frac{f(q-2\eta)}{f(q)}+\prod_{k=1}^n
\frac{z_{u k}(q-2\eta(k-1))}
{\omega_{u k}}\times a_2(q,u)+\\
&&\prod_{k=1}^n \frac{\be_{u k}(\eta,-q)}{\ga_{u k}(q-2\eta(k-1),-)}\times a_3(q,u)\frac{f(q+2\eta)}{f(q)} \ .
\ee
The condition of cancelation is then $K^{(1)}_j=K^{(3)}_j=0 \textrm{ for } 1 \leq j $ and
$K^{(2)}_{lj}=0 \textrm{ for } 1\leq l \leq j$ with the
additional requirement that these three different kinds of condition should in fact lead to the same set
of $n$ nonlinear Bethe equations fixing the $n$ parameters of $\Phi_n$.

Let us first consider the coefficient $K^{(1)}_1$:
\eb
K^{(1)}_1=D_1a_1(u_1)\frac{f(q-2\eta)}{f(q)}+F^{(1)}_1 a_2(q,u_1)
\ee
The condition $K^{(1)}_1=0$ is then equivalent to:
\begin{eqnarray}\label{Bethe_1}
\frac{a_1(u_1)}{a_2(q,u_1)}=\frac{f(q)}{f(q-2\et)}\prod_{k=2}^n
\frac{\tf(u_{k1}+\eta)}{\tf(u_{k1}-\eta)} \times \frac{\tf(q-3\eta)^n}{\tf(q-\eta)^{n-1}\tf(q-2\eta n -\eta)}
\end{eqnarray}
Now it remains to check that the remaining two conditions lead to the same Bethe equations.
The condition
\eb
0=K^{(3)}_1=F^{(2)}_1 a_1(u_1)\frac{f(q)}{f(q+2\eta)}+H_1a_2(q+2\eta)
\ee
yields the same Bethe equation as in \eqref{Bethe_1} thanks to the identity (from the unitarity condition \eqref{unit}):
\eb
\frac{\alpha(\eta,q,u)}{\beta(q,\eta,u)}=-\frac{\alpha(q,\eta,-u)}{\beta(q,\eta,-u)}
\ee
Finally, the cancelation of $K^{(2)}_{12}$ leads also to the same  Bethe equation \eqref{Bethe_1} thanks to
the following identity:
\eb
0&=&\left(\frac{\de_{1u}(-q)}{\ga_{1u}(q,-q)y_{12}(q-2\eta)}+
\frac{z_{1u}(q)\al_{2u}(\eta,q)\omega_{u1}}{\be_{2u}(q,\eta)y_{u1}(q)}\right)\times \tf(q-3\eta)^2+\\
&&\left(\frac{\de_{u1}(q)}{\ga_{u1}(q,-q)y_{12}(q-2\eta)}-\frac{\al_{u1}(q,\eta)}{\ga_{u1}(q,-q)y_{u2}(q-2\eta)}
\right)\times \tf(q-3\eta)^2+\\
&&\frac{1}{y_{u1}(q)}\left(\frac{z_{u1}(q)\al_{u2}(q-2\eta,\eta)}{\be_{u2}(q-2\eta,\eta)}-\frac{\al_{u1}(q,\eta)\al_{12}(q-2\eta,\eta)}
{\be_{u1}(q,\eta)\be_{12}(q-2\eta,\eta)}\right)\times \frac{\tf(u_{21}+\eta)\tf(q-5\eta)\tf(q-\eta)}
{\tf(u_{21}-\eta)}+\\
&&\frac{\al_{u1}(q,\eta)\al_{12}(q+2\eta,\eta)}{\be_{u1}(q,\eta)y_{u1}(q)\be_{12}(q-2\eta,\eta)}\times
\frac{\tf(u_{12}+\eta)\tf(q-5\eta)\tf(q-\eta)}{\tf(u_{12}-\eta)}+\\
&&\frac{\al_{u1}(q,\eta)}{\be_{u1}(\eta,-q)}\left(\frac{z_{u1}(q)}{\omega_{u1}y_{u2}(q)}-
\frac{\al_{u1}(\eta,-q)}{\be_{u1}(q,\eta)y_{12}(q)}\right)\times \frac{\tf(u_{12}+\eta)\tf(q-3\eta)\tf(q-\eta)^2}
{\tf(u_{12}-\eta)\tf(q+\eta)}
\ee

Now we fix $f(q)$ so as to ensure that the Bethe equation (hence its solutions) do not depend on $q$. This
can be achieved by choosing $f(q)=e^{cq}\tf(q-\eta)^n$, where $c$ is an arbitrary constant.

The simultaneous vanishing of $K^{(1)}_1$, $K^{(3)}_1$ and $K^{(2)}_{12}$ then leads to the condition:
\eb
\prod_{k=1}^n\frac{\tf(u_1-z_k-2\eta)}{\tf(u_1-z_k)}=e^{2c\eta}\times
\prod_{k=2}^n \frac{\tf(u_1-u_k-\eta)}{\tf(u_1-u_k+\eta)}
\ee

Once again, the symmetry property of the Bethe vector $\Phi_n$ allows us to derive easily the conditions for
the remaining $u_j$'s by a simple substitution of the spectral parameters. Thus we obtain
 the set of Bethe equations:
\eb
\prod_{k=1}^n\frac{\tf(u_j-z_k-2\eta)}{\tf(u_j-z_k)}=e^{2c\eta}\times
\prod_{\substack{k=1\\k\neq j}}^n \frac{\tf(u_j-u_k-\eta)}{\tf(u_j-u_k+\eta)} \qquad (j=1,\ldots,n)
\ee

Since $f(q)$ is now fixed, we can write the explicit form of the eigenvalues of the transfer
matrix on the module $W=V(z_1)\otimes \ldots \otimes V(z_n)$ as a $q$-independent function of the solutions of the
Bethe equations :
\eb
\Lambda(u,\{u_j\})&=&e^{-2\eta c}\prod_{k=1}^n\frac{\tf(u-z_k-\eta)\tf(u-z_k-2\eta)\tf(u-u_{k}+2\eta)}
{\tf(\eta)\tf(2\eta)\tf(u-u_{k})}+\\
&&+ \prod_{k=1}^n \frac{\tf(u-z_k-\eta)\tf(u-z_k)\tf(u-u_k-2\eta)}{\tf(\eta)\tf(2\eta)\tf(u-u_k)}+\\
&&+e^{2\eta c}\prod_{k=1}^n \frac{\tf(u-z_k)\tf(u-z_k+\eta)\tf(u-u_k-\eta)}{\tf(\eta)\tf(2\eta)\tf(u-u_k+\et)} \, .
\ee

\section{Conclusions}
We showed in this paper that the algebraic Bethe ansatz method can be applied to the elliptic quantum
group $E_{\tau,\eta}(so_3)$. Similarly to Tarasov's implementation of the algebraic Bethe ansatz to the Izergin-
Korepin model, the creation operators for the Bethe vectors  are not simple products of the Lax matrix entries
but are constructed through a recurrence relation. This analogy comes from the fact that
the Izergin-Korepin $R$-matrix, although nondynamical,
has nonzero entries at the same positions as the elliptic dynamical $R$-matrix considered here.
For the simplest highest weight module available we gave the Bethe vectors,
and derived the Bethe equations as well as the eigenvalues of the transfer matrix. Detailed proofs of these
results will be published elsewhere \cite{MaNa2}.
\subsection*{Acknowledgments}

We wish to thank Petr Petrovich Kulish for illuminating discussions.
This work was supported by the project POCI/MAT/58452/2004,
in addition to that Z. Nagy benefited from
the FCT grant SFRH/BPD/25310/2005 and N. Manojlovi\'c acknowledges additional support from
SFRH/BSAB/619/2006. This manuscript was completed at the Centre for Mathematical Science of
the City University, London, we are grateful for their kind hospitality.

\end{document}